\documentclass[twoside]{article}
\usepackage{latexsym}
\usepackage{amsfonts}
\usepackage{graphicx}

\newtheorem{theorem}{Theorem}[section]
\newtheorem{lemma}{Lemma}[section]

\newtheorem{prop}{Proposition}[section]
\newcommand{\qed}{\hfill$\Box$\par\medskip}

\newenvironment{Proof}{\noindent{\sc Proof.}}{\qed}

\def\bhag#1{\noindent
\setcounter{equation}{0}
\section{#1}
}
\def\bfgk#1{{{#1}\kern-5.5pt{#1}}}

\def\HH{{\mathbb H}}

\def\RR{{\mathbb R}}
\def\CC{{\mathbb C}}
\def\ZZ{{\mathbb Z}}
\def\PP{{\mathbb P}}
\def\PPI{{{\rm I}\kern-1pt\Pi}}
\def\SS{{\mathbb S}}

\def\TT{\mathsf T}
\def\a{\alpha}
\def\b #1;{{\bf #1}}
\def\x{{\bf x}}
\def\k{{\bf k}}
\def\y{{\bf y}}

\def\zz{{\mathbf 0}}
\def\e{\epsilon}

\def\O{{\cal O}}

\def\C{{\mathcal C}}
\def\R{{\mathcal R}}

\def\derf#1#2{{#1}^{(#2)}}

\def\esssup{\mathop{\hbox{{\rm ess sup}}}}

\def\be{\begin{equation}}
\def\ee{\end{equation}}
\def\bea{\begin{eqnarray}}
\def\eea{\end{eqnarray}}
\def\eref#1{(\ref{#1})}
\def\disp{\displaystyle}

\def\donchitre#1#2{\vskip 6.5cm\noindent
\parbox[t]{1in}{\special{eps:#1.eps x=6.5cm y=5.5cm}}
\hbox to 7cm{}\parbox[t]{0.0cm}{\special{eps:#2.eps x=6.5cm y=5.5cm}}}

\def\from{From\ }

\title{A construction of linear bounded interpolatory operators on
  the torus}

\author{ S.~Chandrasekaran\thanks{Department of Electrical and
    Computer Engineering, University of California, Santa Barbara,
    Santa Barbara, CA 93106. The research of this author was
    supported, in part, by grants CCF-0515320 and CCF-0830604 from the
    NSF.} \and H.~N.~Mhaskar\thanks{Department of Mathematics, California
    State University, Los Angeles, California, 90032, U.S.A. The
    research of this author was supported, in part, by grant DMS-0605209 and its continuation
    DMS-0908037 from the National Science Foundation and grant
    W911NF-09-1-0465 from the U.S. Army Research Office.}  }

\date{}

\begin{document}
\maketitle
\hspace{1.0cm}
\begin{abstract}
Let $q\ge 1$ be an integer. Given $M$ samples of a smooth function of $q$ variables, $2\pi$--periodic in each variable, we consider the problem of constructing a $q$--variate trigonometric polynomial of spherical degree $\O(M^{1/q})$ which interpolates the given data, remains bounded (independent of $M$) on $[-\pi,\pi]^q$, and converges to the function at an optimal rate on the set where the data becomes dense. We prove that the solution of an appropriate optimization problem leads to such an interpolant. Numerical examples are given to demonstrate that this procedure overcomes the Runge phenomenon when interpolation at equidistant nodes on $[-1,1]$ is constructed, and also provides a respectable approximation for bivariate grid data, which does not become dense on the whole domain. 
\end{abstract}

\bhag{Introduction} Interpolation at equidistant nodes on the unit
interval $[-1,1]$ is a very classical problem. In the first course in
numerical analysis, one learns of the Newton divided difference
algorithm to find such an interpolant, and the corresponding error
formula. The Runge example, $x\mapsto
(x^2+25)^{-1}$, shows that the sequence of these interpolants need not
converge even if the target function is analytic on $[-1,1]$. In
general, Faber's theorem \cite[Theorem~2, p.~27]{natanson} states that
for any interpolation matrix on $[-1,1]$, there exists a continuous
function on $[-1,1]$ such that the corresponding polynomials of
interpolation to this function do not converge.

The situation changes drastically if one allows the degree of the
interpolatory polynomial to be greater than the minimal
required. Thus, the following Theorem~\ref{szabadtheo2} is a simple
consequence of \cite[Theorem~2.7, p.~52]{szabadbook}. For the purpose
of this exposition, we denote the class of all algebraic polynomials
of degree at most $m$ by $\Pi_m$, and define
$\|f\|_{\infty,[-1,1]}:=\sup_{t\in [-1,1]}|f(t)|$. We note that for
$n$ equidistant nodes on $[-1,1]$, the quantity $d_n$ in the following
theorem satisfies $d_n\ge 2/n$.

\begin{theorem}\label{szabadtheo2}
  Let $x_{k,n}=\cos\theta_{k,n}\in [-1,1]$ be an arbitrary system of
  nodes ($0\le\theta_{1,n}<\cdots<\theta_{n,n}\le\pi$) and let
$$
d_n:=\min_{1\le k\le n-1}(\theta_{k+1,n}-\theta_{k,n}).
$$
Then for any $\e>0$, there exist linear polynomial operators $P_n$ on
$C[-1,1]$ with the following properties: {\rm (a)} If $m=\lfloor
\pi(1+\e)/d_n\rfloor$ then $P_n(P)=P$ for all $P\in\Pi_m$, {\rm (b)}
for $f\in C[-1,1]$, $P_n(f)\in \Pi_N$ where $N=(\pi/d_n+1)(1+3\e)$,
{\rm (c)} $P(f,x_{k,n})=f(x_{k,n})$ for $k=1,\cdots,n$, and {\rm (d)}
\be\label{szabadineq} \|f-P_n(f)\|_{\infty,[-1,1]}\le
c\inf_{P\in\Pi_m}\|f-P\|_{\infty,[-1,1]}. \ee
\end{theorem}

In many engineering applications, one has to find a good approximation
to an unknown \emph{multivariate} target function which also
interpolates the function at certain points, sometimes called
landmarks. For example, in the problem of image registration, we are
given a set of locations $x_j \in [-1,1]^2$ in the first image and a
corresponding set of points $y_j \in [-1,1]^2$ in the second
image. The idea is that the location $x_j$ in the first image is the
``same'' as the location $y_j$ in the second image. We then hope to
find a map $g : [-1,1]^2 \to \RR^2$ such that $g(x_j) = y_j$, and
such that $g$ satisfies some smoothness conditions. There are at least
two reasons for insisting on interpolatory approximation in this
situation.  First, the locations might have been chosen at great
costs, including human efforts. Second, if the registration is being
done many times over a sequence of images (for example when we
stitch together video frames to form a large image), then a
non-interpolatory approximation will cause a drift between the first
image and the last image in the sequence.

It is interesting to note that polynomial interpolation in
multivariate setting has a totally different flavor than in the
univariate setting; for example, even if one has exactly as many
points as the dimension of the polynomial space involved, there might
not exist an interpolant from that space. Even if an interpolant
exists, the error bounds for approximation depend heavily on the
geometry of the points. In \cite{approxint}, we proved that an
analogue of Theorem~\ref{szabadtheo2} holds in practically any setting
where the so called direct theorem of approximation holds, provided we
drop the requirement of linearity. In particular, we proved analogous
results in the multivariate setting. However, the results in
\cite{approxint} are not constructive, and do not yield linear
operators.

The purpose of this paper is to develop algorithms to achieve near
best polynomial approximations to smooth multivariate functions, which
satisfy interpolatory constraints. Our constructions will work without
requiring any specific locations for the points where the target
function is evaluated. We refer to such data as \emph{scattered
  data}. We do not require that the data become dense on the whole
cube. In turn, our approximations may not converge on the whole
cube. However, they will converge at the limit points of the data, and we
will estimate the rate of convergence.

To motivate our construction, we revert to the univariate case of
Theorem~\ref{szabadtheo2}. We recall that there is a one to one
correspondence between functions on $[-1,1]$ and even, $2\pi$-periodic
function on $\RR$, given by $f^\circ(\theta)=f(\cos\theta)$. Moreover,
$\|f^\circ\|_{\infty,[-\pi,\pi]}:=\sup_{\theta\in
  [-\pi,\pi]}|f^\circ(\theta)| =\|f\|_{\infty, [-1,1]}$.  Let $r\ge 1$
be an integer, and $f^\circ$ be $r$ times continously differentiable
on $[-\pi,\pi]$. In this discussion, we will write $P_n$ in place of
$P_n(f)$. In view of a theorem of Czipser and Freud \cite{CzipsFreud},
the estimate \eref{szabadineq} implies that
$\|\derf{P_n^\circ}{r}\|_{\infty, [-\pi,\pi]} \le
c\|\derf{f^\circ}{r}\|_{\infty, [-\pi,\pi]}$. Therefore, the
minimization problem ``minimize
$\|\derf{P^\circ}{r}\|_{\infty,[-\pi,\pi]}$ over all $P\in \Pi_N$,
subject to the constraints $P(x_{k,n})=f(x_{k,n})$, $k=1,\cdots,n$''
has a solution $P_n^*$ with the right bounds on the $r$-th derivative
of its periodic version. In view of the Arzela--Ascoli theorem, this
implies that any subsequence of the sequence $\{P_n^*\}$ has a
uniformly convergent subsequence. If $x_0$ is a limit point of a subsequence $\{x_{k,n}\}_{n\in\Lambda}$, then it is not difficult to deduce using the interpolatory conditions
 that $\lim_{n\in\Lambda, n\to\infty}P_n^*(x_0)=f(x_0)$.   In this
paper, we will extend these ideas to the multivariate periodic
setting. Instead of describing the smoothness of the functions in
terms of derivatives, we will consider Sobolev classes. We will also
consider minimization in arbitrary $L^p$ norms; the $L^1$ norm being
of recent interest from the point of view of compressed sensing. Some
technical details, involving a construction of quasi--interpolatory
polynomial operators, are required to prove the rate of convergence of
our constructions.  However, the bulk of the technical details is in
the proof of the feasibility of the optimization problem. We will use
Theorem~2.1 in \cite{approxint} with the appropriate Sobolev spaces,
and will need to prove the analogue of Theorem~3.2 in \cite{approxint}
also with approximation in Sobolev spaces rather than the space of
continuous functions as in that theorem. Our main tool is the
construction of a multivariate analogue of trigonometric polynomial
frames constructed in \cite{trigwav, loctrigwav}.

We state our main results Section~\ref{mainsect}, and illustrate them numerically in Section~\ref{applsect}. The proofs of the results are given in Section~\ref{proofsect}, following some preparation of a technical nature in Section~\ref{prepsect}. At the first reading, it might help to skip this section, referring back to the various statements there on an as needed basis.

We would like to thank \textbf{Karthik Raghuram} for carrying out all the
numerical experiments.

\bhag{Main results}\label{mainsect} 

 In the
sequel, $q\ge 1$ will denote a fixed integer, and we will think of
$2\pi$--periodic functions on $\RR^q$ as functions on $[-\pi,\pi]^q$,
tacitly identified with the $q$ dimensional torus.  Analogous to the univariate case,   any function $f :[-1,1]^q\to\RR$, corresponds uniquely to the $2\pi$--periodic function $f^\circ$ on $\RR^q$ by the correspondence 
$$
f^\circ(\theta_1,\cdots,\theta_q) = f(\cos\theta_1,\cdots,\cos\theta_q).
$$
The symbol
$\|\circ\|$ will denote the Euclidean norm of a vector in $\RR^q$.
Let $\HH^q_n$ denote the class of all trigonometric polynomials in $q$
variables with spherical order at most $n$; i.e.,
$$
\HH^q_n :=\{\sum_{\k\in\ZZ^s,\ \|\k\|\le n}a_\k\exp(i\k\cdot(\circ))\
:\ a_\k\in\CC\}.
$$
Here, we find it convenient to use the same notation even if $n$ is
not an integer. It is not difficult to see that multivariate algebraic polynomials on $[-1,1]^q$ correspond to the trigonometric polynomials of the same order which are symmetric in each of the variables. Therefore, in this paper, we are interested
mainly in the interpolation of multivariate periodic functions; the results can also be applied trivially to the interpolation of functions on $[-1,1]^q$, with suitable smoothness conditions defined in terms of the corresponding periodic function.

If $1\le p\le \infty$, $K\subset [-\pi,\pi]^q$ and $f: K\to\CC$ are
Lebesgue measurable, we write \be\label{normdef}
\|f\|_{p,K}=\left\{\begin{array}{ll}
    \disp\left\{\int_{K}|f(\x)|^pd\x\right\}^{1/p},
    &\mbox{ if $1\le p<\infty$,}\\
    \disp\esssup_{\x\in K}|f(\x)|, &\mbox{ if $p=\infty$.}
\end{array}\right.
\ee The symbol $L^p(K)$ denotes the class of all Lebesgue measurable
functions $f$ for which $\|f\|_{p,K} <\infty$, with the usual
convention that two functions are considered equal if they are equal
almost everywhere. If $K=[-\pi,\pi]^q$, we will omit its mention from
the notations. If $1<p< \infty$, we will write $p':=p/(p-1)$, and
extend this notation to $p=1,\infty$ by setting $1'=\infty$,
$\infty'=1$. If $f\in L^1$, the Fourier coefficients of $f$ are
defined by \be\label{fourcoeffdef} \hat f(\k):=
\frac{1}{(2\pi)^q}\int_{[-\pi,\pi]^q}f(\x)\exp(-i\k\cdot\x)d\x, \qquad
\k\in\ZZ^q.  \ee If $f\in L^p$, then its degree of approximation from
$\HH_n^q$ is defined by
$$
E_{n,p}(f):=\inf_{T\in\HH_n^q}\|f-T\|_p.
$$

If $s\in\RR$, $1\le p\le \infty$, the Sobolev class $W^p_s$ consists
of all $f\in L^p$ for which there exists $\derf{f}{s}\in L^p$ such
that
$$
\widehat{\derf{f}{s}}(\k) =(\|\k\|^2+1)^{s/2}\hat f(\k), \qquad \k\in\ZZ^q.
$$
We define \be\label{sobnormdef} \|f\|_{W^p_s}:=\|\derf{f}{s}\|_p, \ee
and note that $W^p_s$ is a Banach space.  We observe that if $\Delta$
is the Laplacian operator on $\RR^q$, and $s$ is an even, positive
integer, then $\derf{f}{s}=(\Delta+I)^{s/2}f$, where $I$ is the
indentity operator. In particular, in this case, the operator
$f\mapsto \derf{f}{s}$ is a surface derivative operator on the torus
identified with $[-\pi,\pi]^q$. An important property of the spaces
$W^p_s$ is given in the following proposition, which will be proved in
Section~\ref{sobolsect}. Here, and in the rest of this paper, the
symbols $c, c_1,\cdots$ will denote generic positive constants,
depending on such fixed parameters of the problem as $p$, $s$, $q$,
etc. and other quantities explicitly indicated, but their value may
different at different occurrences, even within a single formula. The
notation $A\sim B$ means that $c_1A\le B\le c_2A$.

The following proposition, to be proved in Section~\ref{sobolsect}, gives an integral representation of functions in $W_s^p$.

\begin{prop}\label{sobkernprop} 
  Let $1\le p\le\infty$, $s>q/p$. Then there exists a function $K_s\in
  L^{p'}$ such that \be\label{sobkerndef}
  \widehat{K_s}(\k)=(\|\k\|^2+1)^{-s/2}, \qquad \k\in\ZZ^q.  \ee If
  $f\in W^p_s$, then for almost all $\x\in [-\pi,\pi]^q$,
  \be\label{kernrepresentaion}
  f(\x)=\frac{1}{(2\pi)^q}\int_{[-\pi,\pi]^q}K_s(\x-\y)\derf{f}{s}(\y)d\y
  =\frac{1}{(2\pi)^q}\int_{[-\pi,\pi]^q}K_s(\y-\x)\derf{f}{s}(\y)d\y.
  \ee 
In particular, $f$ is almost everywhere equal to a continuous
  function. Denoting this continuous function again by $f$, we have
  for any $0<s'<s-q/p$, 
\be\label{sobolevbd} E_{2^n,\infty}(f)\le
  c2^{-n(s-q/p)}\|f\|_{W^p_s}, \ \|f\|_\infty\le
  c\|f\|_{W^\infty_{s'}}\le c\|f\|_{W^p_s}.  \ee
\end{prop}

We remark that in the case $p=2$, one can take the following approach
for interpolation of functions in $W^2_s$, $s>q/2$. The
Golomb--Weinberger variation principle \cite{golumb} can be used to
show that the solution of the minimization problem
\be\label{golomb}
\mbox{ minimize } \{\|g\|_{W^2_s}\ :\ g(\y_{j,n})=f(\y_{j,n}), \ j=1,\cdots,M_n\}
\ee
has a solution in the span of $K_{2s}(\circ-\y_{j,n})$, and therefore,
can be found by solving an appropriate system of linear equations. The
stability of this system as well as the error bounds can be estimated
using known techniques from the theory of radial basis functions, for
example, \cite{schaback} (See Theorem~\ref{conditionnotheo} below). However, we are interested in finding
polynomial interpolants for functions in $W^p_s$ for $s>q/p$, without
requiring $p=2$.
 
As customary in the theory of interpolation, let $Y$ be the
interpolation matrix whose $n$-th row $Y_n$ contains $M_n$ vectors
$\{\y_{j,n}\}_{j=1}^{M_n}$. Our theorems will depend upon two
quantities, defined in \eref{meshnormdef} below, that measure the density of these points as well as their
rareness. First, if  $\C\subseteq [-\pi,\pi]^q$ and $\x\in [-\pi,\pi]^q$, we define
$$
\mbox{dist} (\C, \x):=\inf_{\y\in\C}\|\x-\y\|.
$$
Further, if  $K\subseteq
[-\pi,\pi]^q$, we define the mesh norm $\delta(\C,K)$ (respectively,
separation radius $\eta_\C$) of $\C$ by 
\be\label{meshnormdef}
\delta(\C,K):=\sup_{\x\in K}\mbox{dist} (\C, \x), \
\eta_\C:=(1/2)\inf_{\x,\y\in\C,\ \x\not=\y}\|\x-\y\|.  \ee 
We will
simplify our notation, and write $\delta_n(K)$ for $\delta(Y_n,K)$ and
$\eta_n:=\eta(Y_n)$.

Our first theorem, to be proved in Section~\ref{proofsect}, shows the
feasibility of a procedure for finding interpolatory trigonometric
polynomials in $\HH_n^q$.

\begin{theorem}\label{feasibletheo}
  Let $1\le p\le \infty$, $s>q/p$, $Y$ be as above.\\
  {\rm (a)} There exists an integer $N^*$ with $N^*\sim \eta_n^{-1}$
  and a mapping ${\bf P} :W^p_s\to \HH^q_{N^*}$ such that for every
  $f\in W^p_s$, \be\label{sobinterp} {\bf P}(f,\y_{j,n})=f(\y_{j,n}),
  \qquad j=1,\cdots,M_n, \ee and \be\label{sobapprox} \|f-{\bf
    P}(f)\|_{W^p_s}\le c \inf\{ \|f-T\|_{W^p_s}\ :\ T\in
  \HH^q_{N^*}\}.  \ee {\rm (b)} We consider the minimization problem
  \be\label{minproblem} \mbox{ minimize
  }\left\{\frac{1}{{N^*}^q}\sum_{0\le \k\le
      3N^*-1}|\derf{P}{s}(2\pi\k/(3N^*))|^p\right\}^{1/p}, \ee where
  the minimum is over all $P\in\HH_{N^*}^q$, such that
  $P(\y_{j,n})=f(\y_{j,n})$, $j=1,\cdots,M_n$, and an appropriate
  interpretation is understood in the case $p=\infty$.  There exists a
  solution of this problem, $\PP_n^*=\PP_n^*(p,Y_n,f)\in \HH_{N^*}^q$,
  such that $\|\PP_n^*\|_{W^p_s}\le c\|f\|_{W^p_s}$.
\end{theorem}

In practice, it seems that we can take $N^*=4\eta_n^{-1}$ if $p=2$
and $q=1$. We note that the problem \eref{minproblem} has a unique
solution if $p=2$, and the corresponding operator $P_n^*$ is  linear in $f$.

The next theorem, to be proved in Section~\ref{proofsect}, examines the
convergence properties of the sequence $\{\PP_n^*\}$.

\begin{theorem}\label{maintheo}
  Let $1\le p\le \infty$, $s>q/p$, $f\in W^p_s$, $N^*$ and $\PP_n^*$ be
  found as in Theorem~\ref{feasibletheo}.\\
  {\rm (a)} If  $\Lambda$ is a subsequence of positive integers, $\x_0\in [\pi,\pi]^q$, and 
\be\label{limitptcond}
\lim_{n\to\infty\atop n\in \Lambda} \mbox{{\rm dist}} (Y_n,\x_0)=0,
\ee
 then
$$\lim_{n\to\infty \atop n\in \Lambda}\PP_n^*(\x_0)=f(\x_0).$$
{\rm (b)} There exists a constant $\gamma=\gamma(p,q,s)$ (independent of $n$) with the
following property. If $\x_0\in [-\pi,\pi]^q$, and
 $\delta_n([\x_0-\delta,\x_0+\delta]^q)\le
\gamma\delta$, then \be\label{errorbound} \|f-\PP_n^*\|_{\infty,[\x_0-\delta,\x_0+\delta]^q} \le
c\delta^{s-q/p}\|f\|_{W^p_s}.  \ee
\end{theorem}

The proof of Theorem~\ref{feasibletheo} occupies a major part of this
paper. We will use an abstract result from \cite{approxint}, quoted
here as Lemma~\ref{approxintlemma}. To use this result, we need first
to approximate carefully an arbitrary element of the span of
$\{K_s(\circ-\y) : \y\in Y_{N^*}\}$ for a suitable value of $N^*$ by
trigonometric polynomials in $\HH_{N^*}$; indeed, $N^*$ will be
determined so that this approximation works. In turn, this involves an
estimation of the coefficients of this element in terms of the norms
of this element, as well as a good approximation bound on $K_s$.  In
preparation, in Section~\ref{kernsect}, we introduce certain localized
kenels and operators, and prove a number of technical results
concerning these. These enable us after some further preparation to
prove Proposition~\ref{sobkernprop} and study some further properties of the kernel $K_s$ in
Section~\ref{sobolsect}.  The proof of Theorem~\ref{maintheo}(a), as
expected, is a compactness argument. We also need to estimate the
discrete norm used in \eref{minproblem} by the corresponding
continuous norm. The necessary facts are stated in
Lemmas~\ref{sphtrigmzlemma} and \ref{sobcompactlemma}. The proof of
Theorem~\ref{maintheo}(b) is quite simple in the case when $q=1$,
$p=\infty$, and $s$ in an integer: If there are $s$ elements of $Y_n$
in $K=[x_0-\delta,x_0+\delta]$, we take the Lagrange interpolatory
polynomial $L$ for $f$ (and hence, $\PP_{N^*}$) at these points. The
elementary Newton error formula for interpolation yields
$$
\|f-\PP_{N^*}\|_{\infty, K} \le \|f-L\|_{\infty, K}+
\|\PP_{N^*}-L\|_{\infty,K} \le c\delta^{s}\|\derf{f}{s}\|_\infty.
$$
The Newton formula does not hold in the multivariate case, and no
similarly clean estimates are possible independently of the geometry
of the points in question. Therefore, we use a result from
\cite{quasiint2}, quoted as Proposition~\ref{quasiintprop} below, to
construct an analogue of $L$, which is not interpolatory, but utilizes
only the values $f(\y_{j,n})=\PP_{N^*}(\y_{j,n})$. To take care of the
technicalities of noninteger $s$ and $L^p$ norms other than
$p=\infty$, we use the direct and converse theorems of approximation
theory. Although these results are folklore, we could not find them in
the literature in the form which we needed. Therefore, in
Section~\ref{backsect}, we review the results in the form in which we
found them, and reconcile them to our needs. We also describe the
construction of the algebraic polynomial approximation.

\bhag{Numerical experiments}\label{applsect}

In this section, we will  present numerical experiments that demonstrate
the behavior of the method over a wide variety of situations, some of
which do not satisfy the assumptions made in this paper. In the case when $p=2$, the optimization
problem \eref{minproblem} has a numerically effective closed-form solution. In this case, the problem is formulated easier directly in terms of the coefficients of the trigonometric polynomials: 
\be\label{mincoeffproblem}
\mbox{ Find } \arg\min_{\{a_\k\}}\sum_{\|\k\|\le N^*} |a_\k|^2(1+\|\k\|^2)^s,
\ee
subject to the constraints 
$$
\sum_{\|\k\|\le N^*} a_\k\exp(i\k\cdot\y_{j,n}) =f(y_{j,n}), \qquad j=1,\cdots,M_n.
$$
Therefore, we assume in this section that $p=2$, and refer the interpolant resulting as a solution of this problem  as a \emph{minimum Sobolev norm (MNS)} interpolant. In our computations below, we actually consider rectangular sums rather than the spherical sums as in \eref{mincoeffproblem}. The term MNS interpolant will be used for all such minor variations.

We first consider the classical Runge phenomenon by interpolating the
function $f(x)=(1+100x^2)^{-1}$ at equi-spaced points $\{-1+2j/(n+1)\}_{j=1}^n$ on the interval
$[-1,1]$. In Table~\ref{tab:runge} we show the results of our numerical
experiments. The first column of the table shows the number of data
points $n$ used for interpolation.  To avoid any special structure among the points, we chose the values of $n$ as indicated, so as to be essentially (but not exactly) doubling
from step to step.  In all cases, we chose the order of
the interpolatory polynomial to be $2n$, and computed the maximum error by sampling the MSN interpolant at $3n$ equi-spaced points. Columns 2--8 show
the maximum interpolation error with different values of $s$.   The maximum error  decreases more rapidly with increasing $s$, but there are
diminishing returns for higher values of $s$ due to increasing
condition numbers and the concomitant loss of numerical
accuracy. To minimize this loss, we used a special algorithm
  that combines an $LU$ factorization along with the traditional $LQ$
  factorization for solving the minimum norm problem \cite{spie10}.  Our computations show clearly that the interpolants converge; i.e.,
 the Runge phenomenon has disappeared. For comparison
we also show in the last column the approximation error from using
a cubic-spline interpolant.

\begin{table}[h]
\begin{center}
\begin{tabular}{c|cccccc||c}
$n$ & $s=1.5$ & $s=2.5$ & $s=3.5$ & $s=4.5$ & $s=5.5$ & $s=6.5$ & Spline \\
\hline 31 &  3.6212e-03 &  3.1000e-03 &  4.2114e-03  & 5.1510e-03 &  1.0535e-01 &  1.0127e+00 &3.5710e-03\\
61 &  4.5758e-04  & 1.0681e-04 &  4.7994e-05  & 3.0439e-05  & 2.4317e-05 &  2.0644e-04 &6.5167e-04\\
121&   1.7844e-04 &  1.2509e-06 &  8.4610e-08 &  3.2494e-08 &  5.1699e-08 &  2.7949e-07 &4.1035e-05\\
241 &  6.7863e-05 &  2.4084e-07 &  6.8299e-09 &  1.1385e-10 &  2.5882e-09 &  2.7632e-07 &2.3897e-06\\
481 &  2.5210e-05 &  4.5175e-08 &  6.4371e-10 &  5.8196e-11 &  6.5261e-09 &  7.2588e-06 &1.4739e-07\\
961 &  9.2203e-06 &  8.3199e-09 &  6.3179e-11 &  7.5275e-11 &  6.2423e-08 &  7.3983e-04 &9.2473e-09
\end{tabular} 
\end{center}
\caption{Columns 2--7: maximum error for MSN interpolant of Runge's function
  $(1+100x^2)^{-1}$ on $[-1,1]$ with $n$ equi-spaced points and
  different choices of $s$. The order of the interpolant was $2n$
  in all cases. Column 8: maximum error with cubic spline interpolant.}
\label{tab:runge}
\end{table}


Next, we consider a two dimensional interpolation problem on a region
inside the square $[-1,1]\times[-1,1]$. The target function is given in polar
coordinates by
$$
f(r,\theta) :=|r-1/4|^{1/8}|1-r|^{4/5}\sin(r(2\cos\theta+\sin\theta)).
$$ 
 The function is singular on the circles of
radii $0$, $1/4$ and $1$. Furthermore the function does not
satisfy the smoothness conditions of this paper. For the data points, we take those vertices of a square grid which lie in the indicated regions. If $h$ is the length of each side of the squares in this grid, the target polynomial is a bivariate polynomial of coordinatewise degree $\lfloor 2/h\rfloor$.
 In the following tables, $n$ denotes the number of grid points which lie in the region in question, and $m$ is the dimension of the space of interpolatory polynomials. 

In Table~\ref{tab:annu}, we
 compute the maximum error of the interpolant in the annulus $1/2
< r < 3/4$ using approximately $4n$ grid points. The results
are shown in Table~\ref{tab:annu}. This particular annulus is well
removed from the singularities of $f$. Therefore it is pleasing to see
that the MSN interpolant approximates the underlying function very
accurately. We report the maximum error in the annuli $3/4 < r <
19/20$ and $1/4 < r < 3/10$ for the same MSN interpolants in
Tables~\ref{tab:annu2} and~\ref{tab:annu3} respectively. These annuli are
significantly closer to the circles of radii $r=1$ and $r=1/4$ where
the function is singular. Not surprisingly, the error is much larger
here, but still usefully small.

\begin{table}[h]
\begin{center}
\begin{tabular}{rr|cccccccc}
 $n$ & $m$ & $s=1$ & $s=2$ & $s=3$ & $s=4$ & $s=5$ \\ \hline
 100 &  484 & 2.1979e-03 &  1.4012e-03 &  2.0732e-03 &  2.1384e-03 &3.2363e-03 \\
 352 &  1764 & 9.0618e-04 &  1.0095e-03  & 1.5181e-03 &  2.7093e-03 &5.0088e-03 \\
 1280 &  6724 & 2.4002e-04 &  1.4093e-04 &  2.0451e-04 &  2.1890e-04 &3.4952e-04 \\
 4924 &  26244 & 2.3078e-04 &  1.0745e-05 &  2.3017e-05 &  3.4346e-05 &1.4726e-04
\end{tabular} 
\end{center}
\caption{Maximum error of MSN interpolant in region $1/2 < r < 3/4$.}
\label{tab:annu}
\end{table}

\begin{table}[h]
\begin{center}
\begin{tabular}{rr|ccccc}
$n$ & $m$ & $s=1$ & $s=2$ & $s=3$ & $s=4$ & $s=5$ \\ \hline
 100 &  484 & 2.8778e-02 &  4.6506e-03 &  1.2650e-02 &  3.3905e-02 &  8.8795e-02 \\
 352 & 1764  & 3.5019e-03  & 6.8887e-04  & 5.4859e-03  & 2.5252e-02  & 7.2012e-02  \\
 1280 & 6724 & 1.4967e-03  & 3.3349e-04  & 7.7618e-04  & 6.3416e-03  & 2.9345e-02  \\
 4924  & 26244 & 1.4569e-04  & 1.8086e-04  & 4.1918e-04  & 9.2673e-04  & 3.2014e-02  
\end{tabular}
\end{center}
\caption{Maximum error of MSN interpolant in region $3/4 < r < 19/20$.}
\label{tab:annu2}
\end{table}

\begin{table}[h]
\begin{center}
\begin{tabular}{rr|ccccc}
$n$ & $m$ & $s=1$ & $s=2$ & $s=3$ & $s=4$ & $s=5$ \\ \hline
 100 & 484 & 3.1104e-02 &  2.9877e-02 &  2.9187e-02 &  2.9011e-02 &  2.8334e-02 \\
 352 & 1764 & 6.1838e-02 &  5.2437e-02 &  4.8763e-02 &  4.6610e-02 &  4.5226e-02 \\
 1280& 6724 & 1.2503e-02 &  8.0592e-03 &  7.8406e-03 &  8.6545e-03 &  9.2571e-03 \\
 4924& 26244 & 2.0597e-02 &  1.4048e-02 &  1.0375e-02 &  8.6748e-03 &  8.6535e-03
\end{tabular}
\end{center}
\caption{Maximum error of MSN interpolant in region $1/4 < r < 3/10$.}
\label{tab:annu3}
\end{table}

Next, for the same function $f$, we restricted
the samples to the region $\{\{r<1/4\} \cup \{1 < r\}\} \cap \{[-1,1]
\times [-1,1]\}$. Note that this region is essentially made up of 5
pieces. We used $n$ equi-spaced samples in the region and computed the
MSN interpolant with $m$ coefficients for different values of $s$. The
maximum error in the region $r < 1/5$ is reported in
Table~\ref{tab:eannu1}, and in the region $1.1 < r$ in
Table~\ref{tab:eannu2}.

\begin{table}[h]
\begin{center}
\begin{tabular}{rr|ccccc}
$n$ & $m$ & $s=1$ & $s=2$ & $s=3$ & $s=4$ & $s=5$ \\ \hline
 21 &  484 & 7.7266e-03  &  1.4501e-02  &  1.7311e-02   & 1.9202e-02  &  2.0064e-02  \\
 89  &  1764 & 2.1075e-03   & 1.8032e-03   & 2.0722e-03   & 2.1088e-03  &  2.1234e-03  \\
 401  &  6724 & 2.2433e-03    &1.2778e-03  &  7.6602e-04   & 5.5885e-04  &  4.8267e-04  \\
 1637  & 26244 & 8.9254e-04   & 9.2024e-04  &  5.7556e-04  &  2.9398e-04  &  1.5738e-04
\end{tabular}
\end{center}
\caption{Maximum error of MSN interpolant in region $r < 1/5$.}
\label{tab:eannu1}
\end{table}

\begin{table}[h]
\begin{center}
\begin{tabular}{rr|ccccc}
$n$ & $m$ & $s=1$ & $s=2$ & $s=3$ & $s=4$ & $s=5$ \\ \hline
 21 &  484 & 1.9324e-01 &  2.6661e-01 &  5.9105e-01&   1.2583e+00 &  1.5732e+00 \\
 89 & 1764 & 8.6889e-02 &  2.3262e-02  & 4.7827e-02 &  1.5319e-01  & 2.8550e-01  \\
 401&  6724 & 4.4847e-02 &  1.3919e-02 &  9.1716e-02&   3.9371e-01 &  6.6656e+00 \\
 1637&  26244 & 1.6146e-02 &  1.5521e-02&   1.4326e-01&   3.0044e+00 &  7.2167e+01 
\end{tabular}
\end{center}
\caption{Maximum error of MSN interpolant in region $1.1 < r$.}
\label{tab:eannu2}
\end{table}

These experiments show that the proposed scheme can perform well even
on difficult problems, especially in two dimensions where traditional
interpolation schemes require much more work to achieve comparable
accuracy. The proposed method requires special algorithms to execute
efficiently, which will be discussed elsewhere. The ideas presented
here can also be generalized to handle noisy and redundant
observations. These matters will also be reported elsewhere \cite{spie10}.

\bhag{Technical preparation}\label{prepsect}
In this section, we present may technical results which are preparatory to the proof of the main results of Section~\ref{mainsect}. Our proof of Theorem~\ref{feasibletheo} will require Theorem~\ref{sigmaoptheo} and Theorem~\ref{ksapproxtheo}. Subsections~\ref{kernsect} and  \ref{sobolsect} are devoted to the proof of these. In Subsection~\ref{kernsect}, we introduce a localized kernel and the
corresponding operator which will be used throughout this paper, and
prove a number of results regarding these. In particular, we use these results in Subsection~\ref{sobolsect} to prove Proposition~\ref{sobkernprop} and establish a few other facts related to the kernel $K_s$. In Subsection~\ref{backsect}, we review some well known properties of multivariate trigonometric and algebraic polynomial approximation, which will be used in the proof of Theorem~\ref{maintheo}. 

\subsection{Localized kernels}\label{kernsect}
Let $q\ge 1$ be an
integer. For $t>0$, and $h :[0,\infty)\to \RR$, we define formally
\be\label{trigkerndef}
\Psi_t(h,\x):=\sum_{\k\in\ZZ^q}h(\|\k\|/t)\exp(i\k\cdot\x), \qquad
\x\in \RR^q.  \ee We set $\Psi_0(h,\x):=1$ and $\Psi_t(h,\x):=0$ if
$t<0$.

The following theorem summarizes the important localization estimate
for the kernel $\Psi_t$, where we use the notation
$$
{\cal D} f(u)= f'(u)/u.
$$

\begin{theorem}\label{trigkerntheo}
  Let $Q>(q+1)/2$ be an integer, $h : [0,\infty)\to[0,\infty)$ be a
  $Q-1$ times continuously differentiable function supported on $[0,1]$, with an absolutely continuous derivative $\derf{h}{Q-1}$. In addition, we assume that for some
  constants $0<a<b<\infty$, $h(t)=0$ if $t\ge b$, and $h'(t)=0$ if
  $0\le t\le a$.  With $R=(q-1)/2+Q$, we have
  \be\label{trigkernlocalest} 
|\Psi_t(h,\x)|\le c\|{\cal
    D}^Qh\|_{1,[0,\infty)}\frac{t^q}{\min(1, (t\|\x\|)^R)}, \qquad
  \x\in [-\pi,\pi]^q, \ t>0. 
 \ee
 Further, 
\be\label{trigkernsupnorm}
  \max_{\x\in[-\pi,\pi]^q}|\Psi_t(h,\x)| =\Psi_t(h,\zz)\le ct^q\|{\cal
    D}^Qh\|_{1,[0,\infty)}, \qquad t>0, 
\ee
 and for $1\le p\le\infty$,
 \be\label{trigkernlpnorm}
 \|\Psi_t(h,\circ)\|_p \le   ct^{q/p'}\|{\cal D}^Qh\|_{1,[0,\infty)} \qquad t>0.  
\ee
 Here, the
  constants denoted by $c$ may depend upon $a$, $b$, $q$, and $Q$
  only.
\end{theorem}

In order to prove this theorem, we recall that the Bessel function
$J_\a$ can be defined for $\a>-1/2$, $t>0$ by (\cite[Formula~(1.71.6)]{Szego})
\bea\label{besselalphadef}
  J_\a(t)&=&\frac{(t/2)^\a}{\Gamma((2\a+1)/2)\Gamma(1/2)}\int_{-1}^1 e^{itu}(1-u^2)^{\a-1/2}du\nonumber\\
  &=&\frac{(t/2)^\a}{\Gamma((2\a+1)/2)\Gamma(1/2)}\int_{-1}^1 e^{-itu}(1-u^2)^{\a-1/2}du.
\eea
It is customary to define
\be\label{besselhalf}
J_{-1/2}(t)= \frac{(t/2)^{-1/2}}{\Gamma(1/2)}\cos t, \qquad t>0.
\ee
For $f\in L^1(\RR^q)$, we define its inverse Fourier transform by
\be\label{invfourtransdef}
\tilde f(\x)= (2\pi)^{-q}\int_{\RR^q}f(\y)\exp(i\y\cdot\x)d\y, \qquad \x\in\RR^q.
\ee
\begin{lemma}\label{bessellemma}
{\rm (a)} Let $h_0(\x)=h(\|\x\|)$, $\x\in\RR^q$. Then
\be\label{radialfour} \widetilde{h_0}(\x)
=\frac{\|\x\|^{(2-q)/2}}{(2\pi)^{q/2}}\int_0^\infty
h(s)J_{(q-2)/2}(s\|\x\|)s^{q/2}ds.  \ee {\rm (b)} For $\a\ge 1/2$,
\be\label{besselder} \frac{d}{dt} (t^\a J_\a(t)) =t^\a J_{\a-1}(t).
\ee {\rm (c)} We have \be\label{besselbd} |J_\a(t)|\le ct^{-1/2},
\qquad t>0.  \ee
\end{lemma}
\begin{Proof} 
Part (a) is proved, except with a different notation in \cite[Theorem~3.3, p.~155]{steinweiss}. Part (b) is a straightforward consequence of the series expansion for $J_\a$ \cite[Formula~(1.71.1)]{Szego}:
$$
J_\a(t)=\sum_{k=0}^\infty
\frac{(-1)^k}{k!\Gamma(k+\a+1)}(t/2)^{2\a+2k}.  
$$
The estimate \eref{besselbd} follows from \cite[Formula~(7.31.5)]{Szego}.
\end{Proof}

\noindent
\textsc{Proof of Theorem~\ref{trigkerntheo}.} Without loss of
generality, we may assume in this proof that $\|{\cal
  D}^Qh\|_{1,[0,\infty)}=1$. First, we prove
\eref{trigkernsupnorm}. The first equation follows immediately from
the definitions and the fact that $h(t)\ge 0$ for all $t$.  Since $h(\|\k\|/t)=0$ if $\|\k\|\ge bt$, $|h(u)|\le c$ for $u\in\RR$,  and the cardinality of the set $\{\k\in\ZZ^q : \|\k\|\le bt\}$ does not exceed $c_1t^q$, we see from the definition that $0\le \Psi_t(h,0)\le c_2t^q$. This proves the 
last inequality in \eref{trigkernsupnorm}.

In the proof of \eref{trigkernlocalest}, we can assume that
$t\|\x\|\ge 1$.  In this proof only, let $h_0(\x)=h(\|\x\|)$,
$\x\in\RR^q$. In view of the Poisson summation formula
\cite[p.~251]{steinweiss} (our notation is different), we have for
$\x\in [-\pi,\pi]^q$, \be\label{pf2eqn1} \Psi_t(h,\x) = (2\pi)^q
t^q\sum_{\k\in\ZZ^q} \widetilde h_0(t(\x+2\pi\k)).  \ee Let
$\k\in\ZZ^q$, $t(\x+2\pi\k)=\y$, and $\|\y\|=r$. In view of
Lemma~\ref{bessellemma}(a), we have 
\be\label{pf2eqn2}
 \widetilde h_0(\y)=\frac{r^{1-q/2}}{(2\pi)^{q/2}}\int_0^\infty
h(s)J_{(q-2)/2}(sr)s^{q/2}ds.  
\ee 
Let $\a\ge -1/2$. The equation \eref{besselder} used
with $\a+1$ in place of $\a$ shows that 
$$
\int_0^u J_{\a}(rs)s^{\a+1} ds= r^{-\a-2}\int_0^{ru} J_\a(v)v^{\a+1}dv
= r^{-\a-2}(ru)^{\a+1}J_{\a+1}(ru)= \frac{u^{\a+2}J_{\a+1}(u)}{ru}.
$$
Consequently, an integration by parts in \eref{pf2eqn2} yields that
$$
\int_0^\infty h(s)J_{\a}(sr)s^{\a+1}ds=r^{-1}\int_0^\infty {\cal D}
h(u)J_{\a+1}(ru)u^{\a+2}du.
$$
Repeating this $Q$ times, we obtain
$$
\int_0^\infty h(s)J_{\a}(sr)s^{\a+1}ds=r^{-Q}\int_0^\infty {\cal D}^Q
h(u)J_{\a+Q}(ru)u^{\a+Q+1}du.
$$
We recall that $h'(u)=0$ if $u\not\in [a,b]$. Consequently,
$$
\int_0^\infty h(s)J_{\a}(sr)s^{\a+1}ds=r^{-Q}\int_{a}^b {\cal D}^Q
h(u)J_{\a+Q}(ru)u^{\a+Q+1}du.
$$
In view of \eref{besselbd} and the fact that $\a+Q+3/2\ge 0$, we
deduce that 
\be\label{pf2eqn3}
 \left|\int_0^\infty
  h(s)J_{\a}(sr)s^{\a+1}ds\right| \le cr^{-Q-1/2}\|{\cal   D}^Qh\|_{1,[0,\infty)}=cr^{-Q-1/2}, \qquad \a\ge -1/2.  
\ee 
Using $q/2-1$ in place of
$\a$ and substituting the resulting estimate into \eref{pf2eqn2}, we
obtain that \be\label{pf2eqn4} \widetilde h_0(t(\x+2\pi\k))=\widetilde
h_0(\y)\le cr^{1/2-q/2-Q}\|{\cal
  D}^Qh\|_{1,[0,\infty)}=\frac{c}{(t\|\x+2\pi\k\|)^R}.  \ee When
$\k\not=\zz$, we have
$$
\|\x+2\pi\k\|\ge |\x+2\pi\k|_\infty\ge 2\pi|\k|_\infty-\pi\ge
\pi|\k|_\infty\ge \pi\|\k\|/\sqrt{q}.
$$
Since $Q>(q+1)/2$, we have $R>q$, and hence, 
\be\label{pf2eqn5}
\sum_{\k\in\ZZ^q, \ \k\not=\zz} |\widetilde h_0(t(\x+2\pi\k))| \le
ct^{-R}\sum_{\k\in\ZZ^q,\ \|\k\|\ge 1}\frac{1}{\|\k\|^R}\le c(t\|\x\|)^{-R}.  
\ee 
If $\k=\zz$, \eref{pf2eqn4}
yields $|\widetilde h_0(t\x)| \le c(t\|\x\|)^{-R}$. Together with
\eref{pf2eqn5} and \eref{pf2eqn1}, this implies
\eref{trigkernlocalest}.

Since $R>q$, we see from \eref{trigkernlocalest} that
$$
\int_{\|\x\|\ge 1/t} |\Psi_t(h,\x)|d\x \le ct^{q-R}\int_{\|\x\|\ge 1/t}\|\x\|^{-R}d
\x= ct^{q-R}\int_{1/t}^\infty
u^{q-1}u^{-R}du \le c.
$$
Since \eref{trigkernsupnorm} shows that $\int_{\|\x\|\le
  1/t}|\Psi_t(h,\x)|d\x\le c$ as well, we have proved
\eref{trigkernlpnorm} in the case when $p=1$. The estimate
\eref{trigkernlpnorm} in the case $p=\infty$ follows from
\eref{trigkernsupnorm}. The general case is obtained using the
convexity inequality
$$
\|g\|_p \le \|g\|_\infty^{1/p'}\|g\|_1^{1/p}, \qquad g\in L^1\cap
L^\infty,\ 1<p<\infty.
$$
\qed

If $f\in L^1$, we define
\be\label{sigmaopdef}
\sigma_t(h,f,\x):=\frac{1}{(2\pi)^q}\int_{[-\pi,\pi]^q} f(\y)\Psi_t(h,\x-\y)d\y.
\ee
The following theorem summarizes some facts related to this operator.
\begin{theorem}\label{sigmaoptheo}
  Let $h$ satisfy the conditions of Theorem~\ref{trigkerntheo}, $1\le p\le \infty$, $f\in L^p$.\\
{\rm (a)} We have
 \be\label{sigmaopbd}
  \|\sigma_t(h,f)\|_p\le c\|{\cal D}^Qh\|_{1,[0,\infty)}\|f\|_p,
  \qquad t>0.  \ee 
{\rm (b)} In particular, if $h(t)=1$ on $[0,1/2]$ and
  $h(t)=0$ on $[1,\infty)$, then 
\be\label{goodapprox} E_{t,p}(f)\le
  \|f-\sigma_t(h,f)\|_p \le c\left(1+\|{\cal
    D}^Qh\|_{1,[0,\infty)}\right)E_{t/2,p}(f), \qquad t>0.  \ee 
{\rm (c)} If $s>0$ and
  $f\in W^p_s$ then \be\label{favardtheo} E_{t,p}(f)\le
  \frac{c_1}{t^s}\|{\cal D}^Qh\|_{1,[0,\infty)}E_{t,p}(\derf{f}{s}),
  \qquad t>0 \ee and \be\label{bernineq}
  \|\derf{\sigma_t(h,f)}{s}\|_p\le c\|{\cal
    D}^Qh\|_{1,[0,\infty)}t^s\|\sigma_t(h,f)\|_p, \qquad t>0.  \ee 
{\rm (d) (Bernstein inequality)} In
  particular, if $t\ge 0$ and $T\in \HH_t^q$, \be\label{bernineq2}
  \|T\|_{W^p_s}\le ct^s\|T\|_p.  \ee
\end{theorem}

\begin{Proof}\ 
   In view of \eref{trigkernlpnorm},
  $\|\Psi_t(h,\circ)\|_1 \le c\|{\cal D}^Qh\|_{1,[0,\infty)}$. The estimate \eref{sigmaopbd} is now
  clear in the case $p=\infty$, and follows from Fubini's theorem in
  the case when $p=1$. An application of Riesz--Thorin theorem leads
  to the intermediate cases.

  Next, let $h(t)=1$ on $[0,1/2]$ and $h(t)=0$ on $[1,\infty)$. Then
  $\sigma_t(h,T)=T$ for all $T\in \HH_{t/2}^q$. Therefore,
$$
\|f-\sigma_t(h,f)\|_p= \|f-T-\sigma_t(h,f-T)\|_p\le c\left(1+\|{\cal
    D}^Qh\|_{1,[0,\infty)}\right)\|f-T\|_p, \qquad
T\in \HH_{t/2}^q.
$$
This implies \eref{goodapprox}.

Next, let $s\in\RR$, and in this proof only,
$g_t(u)=(h(u)-h(2u))/(u^2+1/t^2)^{s/2}$, $u\in [0,\infty)$, $t\in
(0,\infty)$. Using the fact that $g_t(u)=0$ if $u\in [0,1/4]$, it is not difficult to verify that $g_t$ satisfies
the same conditions as $h$ and $\| {\cal D}^Qg_t\|_{1,[0,\infty)}\le c
\|{\cal D}^Qh\|_{1,[0,\infty)}$ with constant independent of
$t$. Since
$\sigma_t(h,f)-\sigma_{2t}(h,f)=t^{-s}\sigma_t(g_t,\derf{f}{s})$, this
implies that \be\label{pf6eqn1}
\|\sigma_t(h,f)-\sigma_{2t}(h,f)\|_p\le ct^{-s}\|\derf{f}{s}\|_p\|{\cal D}^Qh\|_{1,[0,\infty)}, \qquad t>0, \ s\in\RR.  \ee
Next, let $s>0$. Then \eref{pf6eqn1} leads to 
\begin{eqnarray*}
  E_{t,p}(f)&\le&\|f-\sigma_t(h,f)\|_p=\left\|\sum_{k=0}^\infty (\sigma_{2^kt}(h,f)-\sigma_{2^{k+1}t}(h,f))\right\|_p\\
  &\le& \sum_{k=0}^\infty \|\sigma_{2^kt}(h,f)-\sigma_{2^{k+1}t}(h,f)\|_p \le ct^{-s}\|\derf{f}{s}\|_p\|{\cal D}^Qh\|_{1,[0,\infty)},
\end{eqnarray*}
If $T\in \HH_t$ satisfies $\|\derf{f}{s}-\derf{T}{s}\|_p\le
2E_{t,p}(\derf{f}{s})$, then this implies that
$$
E_{t,p}(f)=E_{t,p}(f-T)\le ct^{-s}\|{\cal
  D}^Qh\|_{1,[0,\infty)}\|\derf{f}{s}-\derf{T}{s}\|_p\le
ct^{-s}E_{t,p}(\derf{f}{s})\|{\cal D}^Qh\|_{1,[0,\infty)}.
$$
This proves \eref{favardtheo}. 

We observe that $\sigma_1(h,f)=\hat f({\bf
  0})=\sigma_1(h,\derf{f}{s})$. We observe also that \eref{pf6eqn1}
holds also for $s<0$. So, if $s>0$, we may apply \eref{pf6eqn1} with
$\derf{f}{s}$ in place of $f$ and $-s$ in place of $s$ to conclude
that
$$
\|\sigma_t(h,\derf{f}{s})-\sigma_{2t}(h,\derf{f}{s})\|_p\le ct^{s}\|f\|_p\|{\cal D}^Qh\|_{1,[0,\infty)}.
$$
 Hence, for $n\ge 1$,
\begin{eqnarray*}
  \|\derf{\sigma_{2^n}}{s}(h,f)\|_p&=&\|\sigma_{2^n}(h,\derf{f}{s}\|_p=\left\|\sigma_1(h,\derf{f}{s})+\sum_{k=0}^{n-1}\left(\sigma_{2^{k+1}}(h, \derf{f}{s})- \sigma_{2^k}(h,\derf{f}{s})\right)\right\|_p\\
  &\le& \|\sigma_1(h,\derf{f}{s})\|_p +\sum_{k=0}^{n-1}\left\|\sigma_{2^{k+1}}(h, \derf{f}{s})- \sigma_{2^k}(h,\derf{f}{s})\right\|_p \\
  &\le& c\|f\|_p\left\{1+\sum_{k=0}^{n-1}2^{ks}\right\}\|{\cal D}^Qh\|_{1,[0,\infty)} \le c2^{ns}\|f\|_p\|{\cal D}^Qh\|_{1,[0,\infty)}.
\end{eqnarray*}
This leads to \eref{bernineq}. The estimate \eref{bernineq2} is
obtained by using \eref{bernineq} with $2t$ in place of $t$ and $T$ in
place of $f$, where we may use a fixed $h$, so that the constant is
independent of the function $h$ used in the rest of the statements of
this theorem.
\end{Proof}

Our next major goal is to prove Theorem~\ref{ksapproxtheo}. In this section, we develop the properties of the kernels $\Psi_n$ which are required in this proof. Let
$\{\y_j\}_{j=1}^M\subset [-\pi,\pi]^q$, $m\ge 1$ be an integer with
\be\label{minsepmconn}
\min_{j\not=k}\|\y_j-\y_k\|\ge 1/m.
\ee
 We note that this implies $M\le
cm^q$. In the sequel, we will assume tacitly that $\{\y_j\}_{j=1}^M$
is one of the members of a sequence of finite subsets of
$[-\pi,\pi]^q$. Thus, $M$ and $m$ are variables, and the constants are
independent of these. If ${\bf a}=\{a_k\}_{k=0}^\infty$ is any sequence of complex numbers, we define
$$
\|{\bf a}\|_{\ell^p} :=\left\{\begin{array}{ll}
\left\{\sum_{k=0}^\infty |a_k|^p\right\}^{1/p}, & \mbox{if $1\le p<\infty$,}\\
\sup_{k\ge 0}|a_k|, & \mbox{if $p=\infty$.}
\end{array}\right.
$$
If ${\bf a}$ is in a Euclidean space $\RR^D$, $\|{\bf a}\|_{\ell^p}:=\|(0,a_1,\cdots,a_D, 0,\cdots)$. 

\begin{prop}\label{coeffprop}
  Let $n\ge 1$ be an integer, $1\le p\le \infty$, ${\bf a}\in\RR^M$,
  $h$, $Q$, $R$ be as in Theorem~\ref{trigkerntheo}, and
  $G(\x) :=\sum_{j=1}^M a_j\Psi_n(h,\x-\y_j)$, $\x\in [-\pi,\pi]^q$. \\
  {\rm (a)} We have \be\label{coeffineq1} \|G\|_p \le
  cn^{q/p'}\left\{1+(m/n)^R\right\}^{1/p'}\|{\cal
    D}^Qh\|_{1,[0,\infty)}\|{\bf a}\|_{\ell^p}.  \ee {\rm (b)} Suppose
  that there exists a compact interval $I\subset (0,1]$ and a constant
  $c_0=c_0(h,I)$ such that $h(t)\ge c_0$ if $t\in I$. Then there
  exists $C_1>0$ depending on $I$, $c_0$, $q$, and $Q$ such that for
  $n\ge C_1m$, \be\label{coeffineq2} c_2n^{-q/p'}\|G\|_p\le \|{\cal
    D}^Qh\|_{1,[0,\infty)}\|{\bf a}\|_{\ell^p} \le
  c_3n^{-q/p'}\|G\|_p.  \ee
\end{prop}

The proof requires a number of preparatory results, some of which we
find of interest in their own right.

\begin{prop}\label{disckernprop}
  Let $h$ satisfy the conditions of Theorem~\ref{trigkerntheo},
  $\{\y_j\}_{j=1}^M\subset [-\pi,\pi]^q$, $m\ge 1$ be an integer with
  $\min_{j\not=k}\|\y_j-\y_k\|\ge 1/m$. For integer $n\ge 1$ and
  $\x\in[-\pi,\pi]^q$, \be\label{disckernloc} \sum_{j, \|\x-\y_j\|\ge
    1/m}|\Psi_n(h,\x-\y_j)| \le cn^q(m/n)^R\|{\cal
    D}^Qh\|_{1,[0,\infty)}.  \ee Hence, \be\label{disckernbd}
  \frac{1}{m^q}\sum_{j=1}^M |\Psi_n(h,\x-\y_j)| \le
  c(n/m)^q\left\{1+(m/n)^R\right\}\|{\cal D}^Qh\|_{1,[0,\infty)}.  \ee
\end{prop}

\begin{Proof}
  Without loss of generality, we may assume that $\|{\cal
    D}^Qh\|_{1,[0,\infty)}=1$. In this proof only, let $\ZZ_k=\{j\ :\
  k/m\le \|\x-\y_j\|\le (k+1)/m\}$, $k=1,2,\cdots$. We note that since
  the minimal separation amongst $\y_j$'s does not exceed $1/m$, there
  are at most $ck^{q-1}$ elements in each $\ZZ_k$. We note that since
  $Q>(q+1)/2$, $R=(q-1)/2+Q>q$.  In view of \eref{trigkernlocalest},
  we have
\begin{eqnarray*}
  \sum_{j, \|\x-\y_j\|\ge 1/m}|\Psi_n(h,\x-\y_j)| &\le &cn^q\sum_{j, \|\x-\y_j\|\ge 1/m}(n\|\x-\y_j\|)^{-R}\\
  &=& cn^{q-R}\sum_{k=1}^\infty \sum_{j\in\ZZ_k}\|\x-\y_j\|^{-R} \le cn^{q-R}m^R\sum_{k=1}^\infty k^{q-1-R}\\
  &\le& cn^q(m/n)^R.
\end{eqnarray*}
This proves \eref{disckernloc}. 

In light of \eref{minsepmconn}, the number of $\y_j$'s with $\|\x-\y_j\|\le 1/m$ is bounded independently of $M$ and $m$.  Hence, \eref{trigkernsupnorm}
implies that
$$
\sum_{j, \|\x-\y_j\|\le 1/m}|\Psi_n(h,\x-\y_j)| \le cn^q.
$$
Together with \eref{disckernloc}, this leads to \eref{disckernbd}.
\end{Proof}

For $f :\{\y_j\}\to\RR$, we will write
$$
\||f\||_p =\left\{\begin{array}{ll}
    \left\{\frac{1}{m^q}\sum_{j=1}^M |f(\y_j)|^p\right\}^{1/p}, &\mbox{ if $1\le p<\infty$,}\\
    \max_{1\le j\le M}|f(\y_j)|, &\mbox{ if $p=\infty$.}
\end{array}\right.
$$

\begin{theorem}\label{mztheo}
Let $1\le p\le \infty$. For any integer $n\ge 1$, and $T\in\HH^q_n$, we have
\be\label{mzineq}
\||T\||_p \le c(n/m)^{q/p}\left\{1+(m/n)^R\right\}^{1/p}\|T\|_p.
\ee
\end{theorem}

\begin{Proof}
  In this proof only, let $h :[0,\infty)\to [0,\infty)$ be a fixed,
  infinitely differentiable function, $h(t)=1$ if $0\le t\le 1/2$,
  $h(t)=0$ if $t\ge 1$, and we choose $Q=q+1$. The constants will
  depend upon this $h$, but $h$ being fixed in this proof, this
  dependence need not be specified.  A comparison of Fourier
  coefficients shows that for $T\in\HH^q_{2n}$,
$$
T(\y)=\frac{1}{(2\pi)^q}\int_{[-\pi,\pi]^q} T(\x)\Psi_{4n}(h,\x-\y)d\y.
$$
In view of \eref{disckernbd}, we obtain
$$
\frac{1}{m^q}\sum_{j=1}^M |T(\y_j)| \le
\|T\|_1\max_{\x\in [-\pi,\pi]^q}\left\{\frac{1}{m^q}\sum_{j=1}^M
  |\Psi_{4n}(h,\x-\y_j)|\right\}\le
c(n/m)^{q}\left\{1+(m/n)^R\right\}\|T\|_1.
$$
If $f\in L^1$, we apply this estimate with $\sigma_{2n}(h,f)$ in place
of $T$, and use Corollary~\ref{sigmaoptheo} (with $p=1$) to deduce
that
$$
\||\sigma_{2n}(h,f)\||_1\le c(n/m)^{q}\left\{1+(m/n)^R\right\}\|f\|_1.
$$
In view of \eref{sigmaopbd} (with $p=\infty$), it is clear that for
$f\in L^\infty$,
$$
\||\sigma_{2n}(h,f)\||_\infty\le \|\sigma_{2n}(h,f)\|_\infty\le c\|f\|_\infty.
$$
An application of Riesz--Thorin interpolation theorem now implies that
for $1\le p\le \infty$, and $f\in L^p$, \be\label{discsigmaopbd}
\||\sigma_{2n}(h,f)\||_p \le
c(n/m)^{q/p}\left\{1+(m/n)^R\right\}^{1/p}\|f\|_p.  \ee If $T\in
\HH^q_n$, then $\sigma_{2n}(h,T)=T$. Therefore, \eref{discsigmaopbd}
implies \eref{mzineq}.
\end{Proof}

Proposition~\ref{matrixinvprop} below is perhaps well known. A proof can be found in \cite[Proposition~6.1]{eignet}. 

\begin{prop}\label{matrixinvprop}
Let $M\ge1$ be an integer, ${\bf A}$ be an $M\times M$ matrix whose $(i,j)$--th entry is $A_{i,j}$. $1\le p\le\infty$, and $\alpha\in [0,1)$. If
\be\label{diagdom}
\sum_{i=1\atop i\not=j}^M |A_{j,i}| \le \alpha |A_{j,j}|,\ \sum_{i=1\atop i\not=j}^M |A_{i,j}| \le \alpha |A_{j,j}|,  \qquad j=1,\cdots,M,
\ee
and $\lambda=\min_{1\le i\le M}|A_{i,i}|>0$, then ${\bf A}$ is invertible, and
\be\label{diagmatrixinvnorm}
\|{\bf A}^{-1}{\bf b}\|_{\ell^p}\le ((1-\alpha)\lambda)^{-1}\|{\bf b}\|_{\ell^p}, \qquad {\bf b}\in\RR^M.
\ee
\end{prop}

We are now in a position to prove Proposition~\ref{coeffprop}.

\noindent
\textsc{Proof of Proposition~\ref{coeffprop}.}
Without loss of generality, we may assume that $\|{\cal D}^Qh\|_{1,[0,\infty)}=1$. In view of \eref{disckernbd}, we have for $\x\in [-\pi,\pi]^q$,
$$
|G(\x)| \le \sum_{j=1}^M |a_j| |\Psi_n(h,\x-\y_j)| \le \|{\bf a}\|_{\ell^\infty}\sum_{j=1}^M |\Psi_n(h,\x-\y_j)|\le  cn^q\left\{1+(m/n)^R\right\}\|{\bf a}\|_{\ell^\infty}.
$$
Thus,
$$
\|G\|_\infty \le cn^q\left\{1+(m/n)^R\right\}\|{\bf a}\|_{\ell^\infty}.
$$
Using \eref{trigkernlpnorm}, we see that
$$
\|G\|_1\le \sum_{j=1}^M |a_j| \|\Psi_n(h,\circ-\y_j)\|_1 \le c\|{\bf a}\|_{\ell^1}.
$$
An application of Riesz--Thorin interpolation theorem with the operator ${\bf a}\mapsto \sum_{j=1}^M a_j\Psi_n(h,\circ-\y_j)$ implies \eref{coeffineq1}. 

Next, if the hypothesis in part (b) is satified, then
\be\label{pf3eqn2}
\Psi_n(h,\zz)\ge \sum_{\k, \ \|\k\|/n\in I} h(\|\k\|/n) \ge cn^q.
\ee
Therefore, \eref{minsepmconn} and \eref{disckernloc} show that for $n\ge C_1m$, $\ell=1,\cdots,M$,
\be\label{pf3eqn1}
\sum_{j=1\atop j\not=\ell}|\Psi_n(h,\y_\ell-\y_j)| \le (1/2)\Psi_n(h,\zz).
\ee
In this proof only, let ${\bf A}$ be the matrix whose $(\ell,j)$-th entry is $\Psi_n(h,\y_\ell-\y_j)$ and ${\bf b}\in\RR^M$ be defined by $b_\ell=G(\y_\ell)$, $\ell=1,\cdots,M$. In view of \eref{pf3eqn1}, \eref{diagdom} is satisfied with $1/2$ in place of $\a$, and in view of \eref{pf3eqn2}, we may choose $\lambda$ to be $cn^q$. Hence, Proposition~\ref{matrixinvprop} implies that ${\bf A}$ is invertible, and 
$$
\|{\bf A}^{-1}{\bf b}\|_{\ell^p}\le cn^{-q}\|{\bf b}\|_{\ell^p}.
$$
Since, ${\bf A}^{-1}{\bf b}={\bf a}$, we have proved that
\be\label{pf3eqn3}
\|{\bf a}\|_{\ell^p} \le cn^{-q}\|{\bf b}\|_{\ell^p}.
\ee
Since $G\in \HH^q_n$, we obtain from Theorem~\ref{mztheo} that
$$
\||G\||_p= m^{-q/p}\|{\bf b}\|_{\ell^p} \le c(n/m)^{q/p}\left\{1+(m/n)^R\right\}^{1/p}\|G\|_p.
$$
Since $n\ge C_1m$, this gives
$$
\|{\bf b}\|_{\ell^p} \le cn^{q/p}\|G\|_p.
$$
Together with \eref{pf3eqn3}, this leads to the second inequality in \eref{coeffineq2}. The first inequality follows from \eref{coeffineq1} and the fact that $n\ge C_1m$.
\qed

\subsection{Sobolev kernel}\label{sobolsect}
Our goal in this section is to prove Proposition~\ref{sobkernprop} and Theorem~\ref{ksapproxtheo}, and establish a few other facts regarding the kernel $K_s$. In particular, we will give in Theorem~\ref{conditionnotheo} an estimate for the norm of the interpolation matrix $K_{2s}(y_j-y_k)$. In the sequel, we assume $Q>(q+1)/2$ is an integer, $h :[0,\infty)\to [0,\infty)$ is a fixed, $Q-1$ times continuously differentiable function with an absolutely continuous derivative $\derf{h}{Q-1}$,  $h(t)=1$ if $0\le t\le 1/2$, $h(t)=0$ if $t\ge 1$, and $h$ is nondecreasing on $[0,\infty)$. We will write $g(t)=h(t)-h(2t)$. Since $h$ is fixed, the dependence of various constants on $h$ need not be indicated.
For $s\in\RR$, we will write
\be\label{tildepsidef}
\tilde\Psi_{n,s}(\x):=\sum_{\k\in\ZZ^q}g(\|\k\|/2^n)(\|\k\|^2+1)^{-s/2}\exp(i\k\cdot\x).
\ee

The following lemma lists some interesting properties of $\tilde\Psi_{n,s}$.
\begin{lemma}\label{tildepsilemma}
Let $s\in\RR$. We have
\be\label{tildepsilocest}
|\tilde\Psi_{n,s}(\x)|\le c\frac{2^{n(q-s)}}{\min(1, (2^n\|\x\|)^R)}, \qquad \x\in [-\pi,\pi]^q.
\ee
Further, 
\be\label{tildepsisupnorm}
\max_{\x\in[-\pi,\pi]^q}|\tilde\Psi_{n,s}(\x)| =\tilde\Psi_{n,s}(\zz)\sim 2^{n(q-s)},
\ee
and for $1\le p\le \infty$,
\be\label{tildepsilpnorm}
\|\tilde\Psi_{n,s}\|_p \le c2^{n(q/p'-s)}.
\ee
\end{lemma}

\begin{Proof}
In this proof only, let $g_n(t)=g(t)/(t^2+1/n^2)^{s/2}$. Then for $\x\in [-\pi,\pi]^q$,
\be\label{pf4eqn1}
\tilde\Psi_{n,s}(\x)= 2^{-ns}\Psi_{2^n}(g_{2^n},\x).
\ee
Each $g_n$ satisfies the conditions of Theorem~\ref{trigkerntheo}, with $a=1/4$, $b=1$. Moreover, $\|{\cal D}^Q g_n\|_{1,[0,\infty)}\le c$. Therefore, all assertions of the lemma, except for the second relation in \eref{tildepsisupnorm}, follow directly from Theorem~\ref{trigkerntheo}. Theorem~\ref{trigkerntheo} also implies that
$\tilde\Psi_{n,s}(\zz)\le c2^{n(q-s)}$. Since $g(1/2)=h(1/2)-h(1)=1$, and $g$ is continuous, there exists a nondegenerate interval $I\subset [1/4,1]$ such that $g(t)\ge 1/2$ if $t\in I$. Hence,
$$
\tilde\Psi_{n,s}(\zz)=\sum_{\k\in\ZZ^q} g(\|\k\|/2^n)(1+\|\k\|^2)^{-s/2}\ge \sum_{\k\in\ZZ^q,\ \|\k\|/2^n\in I}g(\|\k\|/2^n)(1+\|\k\|^2)^{-s/2}\ge c2^{n(q-s)}.
$$
This proves the second relation in \eref{tildepsisupnorm}.
\end{Proof}

\noindent
\textsc{Proof of Proposition~\ref{sobkernprop}.}
 Since $s>q/p$, \eref{tildepsilpnorm} used with $p'$ in place of $p$ shows that
$$
\sum_{n=0}^\infty \|\tilde\Psi_{n,s}\|_{p'} \le c\sum_{n=0}^\infty 2^{n(q/p-s)}<\infty.
$$
So,  the sequence of trigonometric polynomials, defined by
$$
P_N(\x)=1+\sum_{n=0}^N\tilde\Psi_{n,s}(\x)= 1+\sum_{n=0}^N \sum_{\k\in\ZZ^q}g(\|\k\|/2^n)(1+\|\k\|^2)^{-s/2}\exp(i\k\cdot\x)
$$
converges in $L^{p'}$. All the sums in the above expression being finite sums, we obtain for $N\ge 0$,
$$
P_N(\x)=1+\sum_{\k\in\ZZ^q}\sum_{n=0}^N g(\|\k\|/2^n)(1+\|\k\|^2)^{-s/2}\exp(i\k\cdot\x)=\sum_{\k\in\ZZ^q}h(\|\k\|/2^N)(1+\|\k\|^2)^{-s/2}\exp(i\k\cdot\x).
$$
If $\k\in\ZZ^q$, and $2^N\ge 2\|\k\|$, then $h(\|\k\|/2^N)=1$, and $\hat P_N(\k)=(1+\|\k\|^2)^{-s/2}$.
Denoting the $L^{p'}$--limiting function of $P_N$ by $K_s$, it follows that $K_s\in L^{p'}$ and satisfies \eref{sobkerndef}. Moreover, $P_N=\sigma_{2^N}(h,K_s)$, and the bound on $\|\tilde\Psi_{n,s}\|_{p'}$ in \eref{tildepsilpnorm} used with $p'$ in place of $p$ shows that
\be\label{pf7eqn1}
\|K_s-\sigma_{2^N}(h,K_s)\|_{p'}\le c2^{N(q/p-s)}, \qquad N=0,1,2,\cdots.
\ee

Both sides of the first equation in \eref{kernrepresentaion} have the same Fourier coefficients, and hence, they are equal almost everywhere. Similarly, a comparison of Fourier coefficients shows that $K_s(-\x)=K_s(\x)$ for almost all $\x$. This implies the second equation in \eref{kernrepresentaion}. 

For $f\in W^p_s$, a comparison of Fourier coefficients again shows that for integer $m\ge 0$,
$$
\sigma_{2^m}(h,f,\x)=\frac{1}{(2\pi)^q}\int_{[-\pi,\pi]^q}\sigma_{2^m}(h, K_s,\x-\y)\derf{f}{s}(\y)d\y.
$$
So, \eref{pf7eqn1} implies that
$$
\|\sigma_{2^m}(h,f)-\sigma_{2^{m-1}}(h,f)\|_\infty\le \frac{1}{(2\pi)^q}\|\sigma_{2^m}(h,K_s)-\sigma_{2^{m-1}}(h,K_s)\|_{p'}\|\derf{f}{s}\|_p\le c2^{m(q/p-s)}\|\derf{f}{s}\|_p.
$$
Hence, the series $\sigma_1(h,f)+\sum_{m=1}^\infty(\sigma_{2^m}(h,f)-\sigma_{2^{m-1}}(h,f))$ converges uniformly. It is clear that this limit is almost everywhere equal to $f$, and by choosing the continuous representer in the equivalence class of $f$ to be $f$, the limit is $f$. Moreover,
$$
E_{2^n,\infty}(f)\le \|f-\sigma_{2^n}(h,f)\|_\infty \le \sum_{m=n}^\infty \|\sigma_{2^m}(h,f)-\sigma_{2^{m-1}}(h,f)\|_\infty \le 2^{n(q/p-s)}\|\derf{f}{s}\|_p.
$$
This implies the first estimate in \eref{sobolevbd} is now clear. The second set of estimates are proved similarly. \qed

Our proof of Theorem~\ref{feasibletheo} requires the following theorem that describes an approximation of a typical element of the span of $\{K_{s}(\circ-\y_j)\}$. We recall that the solution of the minimization problem \eref{golomb} is in this span (with $2s$ in place of $s$). 

\begin{theorem}\label{ksapproxtheo}
Let $1\le p\le \infty$, $s>q/p$, $\{a_j\}_{j=1}^M\subset\RR$, $G(\x)=\sum_{j=1}^Ma_jK_s(\x-\y_j)$, $\x\in [-\pi,\pi]^q$, and $m\ge 1$ be the smallest integer such that $\min_{\y_j\not=\y_k}\|\y_j-\y_k\|\ge 1/m$. Then there exists an integer $N^*$, independent of $G$, such that $N^*\sim m$ and
\be\label{ksapprox}
\|G-\sigma_{N^*}(h,G)\|_{p'}\le (1/2)\|G\|_{p'}.
\ee
\end{theorem}

\begin{Proof}\ 
As in the proof of Lemma~\ref{tildepsilemma}, in this proof only, we write $g_n(t)=g(t)/(t+1/n)^s$. Then 
each $g_n$ satisfies the conditions of Theorem~\ref{trigkerntheo}, with $a=1/4$, $b=1$. Moreover, $\|{\cal D}^Q g_n\|_{1,[0,\infty)}\sim 1$, and \eref{pf4eqn1} holds. In this proof only, let 
$$
G_n(\x):=\sum_{j=1}^M a_j\tilde\Psi_{n,s}(\x-\y_j)=\sigma_{2^n}(h, G,\x)-\sigma_{2^{n-1}}(h, G, \x), \qquad \x\in [-\pi,\pi]^q.
$$ 
Then \eref{sigmaopbd} implies that $\|G_n\|_{p'}\le c\|G\|_{p'}$.  Moreover, the proof of Proposition~\ref{sobkernprop} shows that
\be\label{pf4eqn3}
G(\x)-\sigma_{2^N}(h,G,\x)=\sum_{n=N}^\infty G_n(\x),
\ee
with convergence in the sense of $L^{p'}$. 

In view of \eref{pf4eqn1}, \eref{coeffineq2} applied with $\Psi_{2^n}(g_n)$  yields that for $n\ge \log_2(C_1m)$,
\be\label{pf4eqn2}
c_22^{n(s-q/p)}\|G_n\|_{p'}\le \|{\bf a}\|_{\ell^{p'}} \le c2^{n(s-q/p)}\|G_n\|_{p'}.
\ee
We now choose $L$ so that $2^L$ is the smallest power of $2$ exceeding $C_1m$.  Then the second inequality in \eref{pf4eqn2}, used with $L$ in place of $n$,   gives 
\be\label{pf4eqn4}
 \|{\bf a}\|_{\ell^{p'}} \le cm^{(s-q/p)}\|G_L\|_{p'} \le cm^{(s-q/p)}\|G\|_{p'}.
\ee
\from \eref{pf4eqn3}, \eref{pf4eqn2}, and \eref{pf4eqn4}, we conclude that for $2^N\ge C_1m$,
$$
\|G-\sigma_{2^N}(h,G)\|_{p'}\le \sum_{n=N}^\infty\|G_n\|_{p'}\le c\|{\bf a}\|_{\ell^{p'}}\sum_{n=N}^\infty 2^{-n(s-q/p)}\le c(m2^{-N})^{(s-q/p)}\|G\|_{p'}.
$$
We now choose $N$ so that $2^N\sim m$ and the last term above is at most $(1/2)\|G\|_{p'}$, and set $N^*=2^N$.
\end{Proof}

We note a consequence of the proof, which might be of independent interest in view of the fact that the interpolant which yields the minimal Sobolev norm amongst all interpolants is in the span of $\{K_{2s}(\circ-\y_k)\}_{k=1}^M$. The following theorem gives the norm of the inverse of the interpolation matrix $(K_{2s}(y_j-y_k))$  in terms of the minimal separation (equivalently, $m$).
\begin{theorem}\label{conditionnotheo}
Let $s>q/2$, and ${\cal I}$ be the $M\times M$ matrix whose $(j,k)$-th entry is $K_{2s}(y_j-y_k)$, where the points $\y_j$ satisfy \eref{minsepmconn}. Then ${\cal I}$ is positive definite, and
\be\label{invinterpest}
\|{\cal I}^{-1}\|\le cm^{s-q/2}.
\ee
\end{theorem}

\begin{Proof}\ 
We observe that a comparison of Fourier coefficients shows that
$$
K_{2s}(\y_j-\y_k) = \frac{1}{(2\pi)^q}\int_{[-\pi,\pi]^q}K_s(\x-\y_j)K_s(\x-\y_k)d\y.
$$
Let ${\bf a}\in \RR^M$, and $G =\sum_{j=1}^M a_j K_s(\circ-\y_j)$. Then the above identity leads to
\be\label{pf8eqn1}
\sum_{k,j=1}^M a_ja_kK_{2s}(\y_j-\y_k) =\frac{1}{(2\pi)^q}\int_{[-\pi,\pi]^q}\left|\sum_{j=1}^M a_jK_s(\x-\y_j)\right|^2d\y =\|G\|_2^2.
\ee
The estimate \eref{pf4eqn4} used with $p=2$ now shows that
$$
\sum_{k,j=1}^M a_ja_kK_{2s}(\y_j-\y_k)\ge cm^{-(s-q/2)}\|{\bf a}\|_{\ell^2}.
$$
Thus, ${\cal I}$ is a positive definite matrix. In view of the Raleigh-Ritz theorem \cite[Theorem~4.2.2, p.~176]{horn}, the lowest eigenvalue of this matrix is at least $cm^{-(s-q/2)}$. This implies \eref{invinterpest}.
\end{Proof}

Although not strictly a property of the kernels $K_s$, we find it convenient to record the following lemma, which will be needed in our proof of  Theorem~\ref{maintheo}. This lemma  is proved in  much greater generality in \cite[Theorem~3.2, Chapter~15]{Lorentz}.
\begin{lemma}\label{sobcompactlemma}
Let $1\le p\le \infty$, $s'>0$. Then for any $c>0$, the set $B_{c,s',p}:=\{f\in L^p\ : \ \sup_{n\ge 1} 2^{ns'}E_{n,p}(f)\le c\}$ is compact in $L^p$.
\end{lemma} 

\subsection{Background on approximation theory}\label{backsect}
The proof of Theorem~\ref{maintheo} depends upon a number of facts from classical approximation theory, as well as some recent developments. In this section, we review the necessary facts.

First, for integer $r\ge 1$, the modulus of smoothness $\omega_r(f,\delta)$ of a $2\pi$--periodic univariate function $f\in L^p([-\pi,\pi])$ is defined by first defining the forward difference operator
$$
\Delta^r_t f(x)= \sum_{k=0}^r {r\choose k}(-1)^k f(x+kt),
$$
and setting
$$
\omega_r(f,\delta) :=\max_{|t|\le \delta}\|\Delta^r_t f\|_p.
$$
If $f\in L^p([-\pi,\pi]^q)$ and ${\bf r}=(r_1,\cdots,r_q)\ge 0$, ${\bf r}\not=(0,\cdots,0)$,  is a multi--integer, then the modulus of smoothness is defined in \cite[Section~3.4.34]{timan} by 
$$
\omega_{{\bf r}}(f, {\bf h}) = \max_{|t_1|\le h_1, \cdots, |t_q|\le h_q}\left\|\left(\prod_{k=1}^q\Delta_{t_k,k}^{r_k}\right) f\right\|_p,
$$
where the notation $\Delta_{t_k,k}^{r_k}$ means that the operator $\Delta_{t_k}^{r_k}$ is applied to the $k$-th variable in the argument of $f$ and $\Delta_{t_k,k}^0f$ means that no difference is taken with respect to th $k$--th variable. We will write ${\bf e}_k$ to denote the vector in $\RR^q$ with $k$-th coordinate equal to $1$ and the remaining coordinates equal to $0$.


For an integer $n\ge 0$, the class of (rectangular) trigonometric polynomials of order at most $n$ is defined by 
$$
\HH_n^{q,R}:=\left\{\sum_{\k\in\ZZ^q,\ |k_\ell|\le n, \ \ell=1,\cdots,q}a_\k\exp(i\k\cdot(\circ))\ :\ a_\k\in\CC\right\}.
$$
For $f\in L^p$, the degree of approximation from $\HH_n^{q,R}$ is defined by
$$
E_{n,p;R}(f)=\inf_{P\in\HH_n^{q,R}}\|f-P\|_p.
$$
Let $r\ge 0$ be an integer, $\a\in (0,1]$.It is proved in \cite[Section~3.6.4]{timan} that for $f\in L^p$, the relation 
\be\label{approxbdtiman}
\sup_{n\ge 1}n^{r+\a}E_{n,p;R}(f) <\infty
\ee
holds if and only if $f$ has almost everywhere defined partial derivatives $D_k^{r}f$ satisfying 
\be\label{smoothtiman}
\max_{1\le k\le q}\sup_{\delta>0}\delta^{-\a}\omega_{2{\bf e}_k}(D_k^{r}f, \delta{\bf e}_k)<\infty.
\ee
(The formulation in \cite{timan} is not quite precise. However, the version which we have stated can be obtained using the same ideas as in \cite{timan}. See \cite[4.5.6]{trigub} for an analogous statement in a slightly different context.)  
Since $\HH_n^q\subseteq \HH_n^{q,R}\subseteq \HH_{\sqrt{q}n}^q$, we have
$$
E_{\sqrt{q}n,p}(f) \le E_{n,p;R}(f)\le E_{n,p}(f).
$$
Thus, \eref{approxbdtiman} is equivalent to
$$
\sup_{n\ge1}n^{r+\alpha}E_{n,p}(f) <\infty.
$$
We summarize these observations in the following proposition.
  
\begin{prop}\label{timaneqprop}
Let $1\le p\le\infty$, and $f\in L^p$, $r\ge 0$ be integer, and $\a\in (0,1]$. Then
\be\label{approxbdsphtrig}
\sup_{n\ge 1}n^{r+\a}E_{n,p}(f)<\infty
\ee
if and only if $f$ has almost everywhere defined partial derivatives $D_k^{r}f$ satisfying \eref{smoothtiman}.
\end{prop}

Next, we recall some results from the theory of algebraic polynomial approximation. Let $\Pi_r^q$ denote the set of all algebraic polynomials of coordinatewise degree at most $r$. We wish to construct an approximation to a continuous function $f$ on $[-1,1]^q$, defined analogously to \eref{normdef}, based on an arbitrary data set  $\C\subset [-1,1]^q$. The mesh norm $\delta_\C :=\delta(\C, [-1,1]^q)$ is defined analogously to \eref{meshnormdef}. We divide $[-1,1]^q$ into equal subcubes of side $2\delta_\C$; the set of these subcubes will be denoted, in this part of the discussion only, by ${\mathcal R}_\C$. Each of the subcubes has at least one point of $\C$. We form a subset $\C_1\subset \C$ by choosing exactly one point $\xi$ in each $R_\xi\in {\mathcal R}_\C$. Then it is clear that $\delta_{\C_1}\sim \delta_\C$. In the following discussion, the points in $\C\setminus\C_1$ do not play any role, and accordingly, we rename $\C_1$ to be $\C$. The following proposition follows from \cite[Theorem~3.1]{quasiint2}, by taking the functional $P\to\int_{[-1,1]^q}P(\x)d\x$ in place of $\gamma$ in that theorem.

\begin{prop}\label{quasiintprop}
Let  $\C$, $\R_\C$ be as above, $r\ge 1$ be an integer. There exists a constant
$\gamma :=\gamma(r,q)$ with the following property. If $\delta_\C \le \gamma$, then 
\be\label{polymzineq}
\sum_{\xi\in\C}\hbox{{\rm vol}}_q(R_\xi)|P(\xi)| \sim \int_{[-1,1]^q}|P(\x)|d\x, \qquad P\in\Pi_{2r}^q.
\ee
Further,
there exist real numbers $\{a_\xi\ :\ \xi\in\C\}$, such that 
\be\label{weightbd}
|a_\xi|\le c\thinspace\hbox{{\rm vol}}_q(R_\xi)\le c\delta_\C^q, \qquad \xi\in\C,
\ee
and
\be\label{polyquad}
\sum_{\xi\in\C}a_\xi P(\xi) = \int_{[-1,1]^q}P(\x)d\x, \qquad P\in\Pi_{2r}^q.
\ee
\end{prop}

In this section only, let $p_k$ denote the orthonormalized Chebyshev
polynomial of degree $k$, with positive leading coefficient. We define
$$
v_r(x,y)=\sum_{k=0}^{2r} p_k(x)p_k(y), \qquad x,y\in[-1,1],\ r=1,2,\cdots.
$$
and  extend this definition by writing
$$
v_r(\x,\y)=\prod_{\ell=1}^q v_r(x_\ell,y_\ell), \qquad \x,\y\in [-1,1]^q,
$$
If $f :[-1,1]^q\to\RR$ is continuous, $\C$ and $\{a_\xi\}$ are as in
Proposition~\ref{quasiintprop}, we define \be\label{vrdef} V_r(f,\x)=
\sum_{\xi\in\C} a_\xi f(\xi)v_r(\x,\xi).  \ee Using \eref{weightbd},
\eref{polymzineq}, we conclude that \be\label{vmopbd}
\|V_r(f)\|_{\infty,[-1,1]^q} \le c(r)\|f\|_{\infty,[-1,1]^q}.  \ee In
view of \eref{polyquad}, $V_r(P)=P$ for all $P\in\Pi_r^q$. So,
choosing $P^*\in\Pi_r^q$ with $\|f-P^*\|_{\infty,[-1,1]^q} \le
2\inf_{P\in\Pi_r^q}\|f-P\|_{\infty,[-1,1]^q}$, \eref{vmopbd} yields
\be\label{polyapproxbd}
\|f-V_r(f)\|_{\infty,[-1,1]^q}=\|f-P^*-V_r(f-P^*)\|_{\infty,[-1,1]^q}\le
c\|f-P^*\|_{\infty,[-1,1]^q} \le
c\inf_{P\in\Pi_r^q}\|f-P\|_{\infty,[-1,1]^q}.  \ee Using the direct
theorem of approximation theory \cite[Section~5.3.1]{timan}, we
conclude that if $f$ has continuous partial derivatives of order up to
$r$, then \be\label{finalpolyapproxbd}
\|f-V_r(f)\|_{\infty,[-1,1]^q}\le c(r)\sum_{k=1}^q\omega_{2{\bf
    e}_k}(D_k^{r}f, {\bf e}_k/r), \ee where the modulus of smoothness
is defined analogously to \eref{smoothtiman}, except that the maximum
is taken for only those values of $t_1,\cdots, t_k$ which don't take
the argument out of the cube in question. We note again that the
operator $V_r$ is determined entirely by the values
$\{f(\xi)\}_{\xi\in\C}$.

We end this section by recording another observation, establishing a connection between the discrete norm used in the statement of the minimization problem \eref{minproblem} and the continuous $L^p$ norm, which will be needed in the proof of Theorem~\ref{feasibletheo}(b).

\begin{lemma}\label{sphtrigmzlemma}
For integer $n\ge 1$, $1\le p\le\infty$, and $T\in \HH_n^q$, we have
\be\label{sphtrigmzineq}
\left\{\frac{1}{n^q}\sum_{0\le \k\le 3n-1}|T(2\pi\k/(3n))|^p\right\}^{1/p}\sim \|T\|_p.
\ee
\end{lemma}
\begin{Proof}\ 
When $q=1$, \eref{sphtrigmzineq} is the classical Marcinkiewicz--Zygmund inequality \cite[Chapter~X, Theorems~7.5, 7.28]{zygmund}. If $T\in \HH_n^q$, then $T\in \HH_{(n,\cdots,n)}$. So, in the case when $q>1$, one obtains \eref{sphtrigmzineq} by applying its univariate version to each of the variables separately.
\end{Proof}

\bhag{Proofs of the main results in Section~\ref{mainsect}.}\label{proofsect}

Our proof of Theorem~\ref{feasibletheo}(a) relies upon the next lemma,  proved in \cite[Theorem~2.1]{approxint}.

\begin{lemma}\label{approxintlemma}
Let $X$ be a normed
linear space,  $V\subset X$ be a finite dimensional subspace of $X$, $X^*$ be the dual space of $X$, $\{x^*_j\}_{j=1}^M\subset X^*$, and  $Z_*$ be the span of $\{x^*_j\}_{j=1}^M$. Suppose that the restriction map $S:
z^*\in Z_*\mapsto z^*|_{V}$ is injective, and the operator norm $\|S^{-1}\|\le \kappa$ for some $\kappa>0$. Then for every $f\in X$ and $\kappa_1>\kappa$, 
there exists $\TT(f)\in V$ such that 
\be\label{abstract_interp}
z^*(\TT(f))=z^*(f) \qquad \mbox{\rm for every } z^*\in Z_*,
\ee
and
\be\label{abstract_approx}
\|f-\TT(f)\|_X\le (1+\kappa_1)\inf_{v\in V}\|f-v\|_X.
\ee
\end{lemma}

We will use this lemma with $W^p_s$ in place of $X$, $\HH^q_{N^*}$ in place of $V$ for a suitable $N^*$, and point evaluation functionals in place of $x^*_j$'s.

\noindent
\textsc{Proof of Theorem~\ref{feasibletheo}.}
The proof of this theorem is similar to that of \cite[Theorem~3.1]{approxint}, except that the details are much more complicated, requiring the use of Theorem~\ref{sigmaoptheo} and Theorem~\ref{ksapproxtheo}. In this proof only, let $X=W^p_s$, $x_j^*(f)=f(\y_j)$, $j=1,\cdots,M$. Since $s>q/p$, Proposition~\ref{sobkernprop} implies that $x_j^*\in X^*$, $j=1,\cdots,M$. Let $\{a_j\}_{j=1}^M \subset \RR$ and $z^*=\sum_{j=1}^M a_jx_j^*$. Let $N^*$ be as in Theorem~\ref{ksapproxtheo}, and $V=\HH^q_{N^*}$. To estimate $\|S^{-1}\|$ for the operator $S$ as in Lemma~\ref{approxintlemma}, we need to find $T\in \HH_{N^*}^q$ for a suitable $N^*$, and estimate $|z^*(T)|/\|T\|_{W_s^p}$ from below. Let $f^*$ be chosen so that $\|z^*\|_{X^*}\le (4/3)|z^*(f^*)|$ and $\|f^*\|_{W^p_s}=1$. We will prove that $\sigma_{N^*}(h,f^*)\in\HH_{N^*}^q$ (cf. \eref{sigmaopdef}) satisfies 
\be\label{pf5eqn4}
\sup_{T\in V}\frac{|z^*(T)|}{\|T\|_{W^p_s}} \ge \frac{|z^*(\sigma_{N^*}(h,f^*))|}{\|\sigma_{N^*}(h,f^*)\|_{W^p_s}} \ge c\|z^*\|_{X^*}.
\ee
The part (a) of the theorem will then follow from Lemma~\ref{approxintlemma}.

We start by observing that for $f\in W^p_s$, \eref{kernrepresentaion} shows that
$$
z^*(f)=\frac{1}{(2\pi)^q}\int_{[-\pi,\pi]^q}\left\{\sum_{j=1}^M a_jK_s(\y-\y_j)\right\}\derf{f}{s}(\y)d\y=\frac{1}{(2\pi)^q}\int_{[-\pi,\pi]^q}G(\y)\derf{f}{s}(\y)d\y,
$$
where $G$ is defined as in Theorem~\ref{ksapproxtheo}. In light of the duality principle and the definition \eref{sobnormdef}, we see that
\be\label{pf5eqn1}
\|z^*\|_{X^*}=\sup\{ |z^*(f)|\ : \ \|f\|_{W^p_s}=1\}=\|G\|_{p'}.
\ee
 Further a comparison of Fourier coefficients implies that for any integer $n$,
\begin{eqnarray*}
z^*(\sigma_{n}(h,f))&=& \frac{1}{(2\pi)^q}\int_{[-\pi,\pi]^q}G(\y)\derf{(\sigma_{n}(h,f))}{s}(\y)d\y=\frac{1}{(2\pi)^q}\int_{[-\pi,\pi]^q}G(\y)\sigma_{n}(h,\derf{f}{s},\y)d\y\\
&=&\frac{1}{(2\pi)^q}\int_{[-\pi,\pi]^q}\sigma_{n}(h,G,\y)\derf{f}{s}(\y)d\y
\end{eqnarray*}
 Then 
\begin{eqnarray*}
\lefteqn{|z^*(f^*)-z^*(\sigma_{N^*}(h,f^*))|=\left|\frac{1}{(2\pi)^q}\int_{[-\pi,\pi]^q}\left(G(\y)-\sigma_{N^*}(h,G,\y)\right)\derf{f^*}{s}(\y)d\y\right| }\\
&\le& \|G-\sigma_{N^*}(h,G)\|_{p'}\le (1/2)\|G\|_{p'}=(1/2)\|z^*\|_{X^*}\le (2/3)|z^*(f^*)|.
\end{eqnarray*}
Therefore, 
\be\label{pf5eqn2}
|z^*(\sigma_{N^*}(h,f^*))|\ge (1/3)|z^*(f^*)|\ge (1/4)\|z^*\|_{X^*}.
\ee
Moreover, \eref{sigmaopbd} implies that
\be\label{pf5eqn3}
\|\sigma_{N^*}(h,f^*)\|_{W^p_s} =\|\derf{(\sigma_{N^*}(h,f^*))}{s}\|_p =\|\sigma_{N^*}(h,\derf{f^*}{s})\|_p\le c\|\derf{f^*}{s}\|_p =c.
\ee
The estimate \eref{pf5eqn4} follows from \eref{pf5eqn2} and \eref{pf5eqn3}.

We note that necessarily, $\|{\bf P}(f)\|_{W^p_s}\le \|f\|_{W^p_s}$. Therefore, part (b) is  a simple consequence of  Lemma~\ref{sphtrigmzlemma}. 
\qed

\noindent
\textsc{Proof of Theorem~\ref{maintheo}.} 

To prove part (a), we observe that in view of \eref{sobolevbd} and the
fact that $\|\PP_n^*\|_{W^p_s}\le c\|f\|_{W^p_s}$ for all $n$, the
sequence $\{\PP_n^*\}\subset B_{c,s-q/p,\infty}$ for a suitable
constant $c$.  Let $\Lambda_1$ be any subsequence of $\Lambda$. Then Lemma~\ref{sobcompactlemma} shows that the sequence
$\{\PP_n^*\}_{n\in\Lambda_1}$ has a subsequence
$\{\PP_n^*\}_{n\in\Lambda_2}$, which converges uniformly. Let $P$ be
the limit of this subsequence. We will show that if \eref{limitptcond} is satisfied, then $P(\x_0)=f(\x_0)$. Let $\e>0$ be arbitrary. Since
$P$ and $f$ are continuous on $[-\pi,\pi]^q$, there is
$\tilde\delta>0$ such that
$$
|f(\x)-f(\y)|\le \e/3, \ |P(\x)-P(\y)|\le \e/3, \qquad \mbox{ for all
} \x,\y\in [-\pi,\pi]^q,\ \|\x-\y\|\le \tilde\delta.
$$
Further, there exists $N$ so that $n\ge N$, $n\in \Lambda_2$ imply
that $\|P-\PP_n^*\|_\infty\le \e/3$. In view of \eref{limitptcond}, there exists $n\in\Lambda_2$, $n\ge N$ such that some point
$\y_{j,n}\in Y$ satisfies $\|\y_{j,n}-\x_0\|\le \tilde\delta$. Then
$f(\y_{j,n})=\PP_n^*(\y_{j,n})$, and we have
\begin{eqnarray*}
  |f(\x_0)-P(\x_0)|&\le& |f(\x)-f(\y_{j,n})| +|f(\y_{j,n})-\PP_n^*(\y_{j,n})|+|\PP_n^*(\y_{j,n})-P(\y_{j,n})|+|P(\y_{j,n})-P(\x_0)|\\
  &\le& \e/3+0+\|\PP_n^*-P\|_\infty+\e/3\le \e.
\end{eqnarray*}
Since this is true for every subsequential limit of $\PP_n^*$, $n\in\Lambda$, this
proves part (a).

To prove part (b), let $r$ be an integer and $\a\in(0,1]$ be chosen so
that $s-q/p=r+\a$. Since $\PP_n^*, f\in B_{c,s-q/p,\infty}$,
Proposition~\ref{timaneqprop} implies that they both have $r$
derivatives satisfying \eref{smoothtiman}. In this proof only, let
$\PP(\y)=\PP_n^*(\x_0+\delta\y)$, $\tilde f(\y)=f(\x_0+\delta\y)$,
$\y\in [-1,1]^q$. Then the assumptions of part (b) ensure that we can
construct the operator $V_r$ as in \eref{vrdef} based on $Y_n\cap K$ in place of $\C$. In this proof only, if $k\in\{1,\cdots,q\}$, $F=D_k^rf$, then $D_k^r\tilde f(\y)=\delta^rF(\x_0+\delta\y)$. Further,
$$
\Delta_{1/r,k}^2 (D_k^r\tilde f)(\y)=(\Delta_{\delta/r,k}^2F)(\x_0+\delta \y).
$$
Hence,
$$
\|\Delta_{1/r,k}^2 (D_k^r\tilde f)\|_{\infty,[-1,1]^q} = \|\Delta_{\delta/r,k}^2F(\x_0+\circ)\|_{\infty,[-\delta,\delta]^q} \le \omega_2(D_k^rf,(\delta/r) {\bf e}_k) \le c\delta^\a.
$$
A similar estimate holdes also for $\PP$ in place of $\tilde f$.
Using
\eref{finalpolyapproxbd} and the fact that $r+\a=s-q/p$, we deduce that
$$
\|\tilde f-V_r(\tilde f)\|_{\infty,[-1,1]^q}\le c\delta^{s-q/p}, \
\|\PP-V_r(\PP)\|_{\infty,K}\le c\delta^{s-q/p}.
$$
We now observe that $V_r(\tilde f)=V_r(\PP)$, and hence, the above
inequalities imply that $\|\tilde f-\PP\|_{\infty, [-1,1]^q}\le
c\delta^{s-q/p}$. Scaling back to the original scale, we obtain
\eref{errorbound}. \qed


\end{document}